# The Distribution of Mixing Times in Markov Chains

Jeffrey J. Hunter

*School of Computing & Mathematical Sciences, Auckland University of Technology, Auckland, New Zealand*

December 2010

**Abstract**

The distribution of the "mixing time" or the "time to stationarity" in a discrete time irreducible Markov chain, starting in state *i*, can be defined as the number of trials to reach a state sampled from the stationary distribution of the Markov chain. Expressions for the probability generating function, and hence the probability distribution of the mixing time starting in state *i* are derived and special cases explored. This extends the results of the author regarding the expected time to mixing [J.J. Hunter, Mixing times with applications to perturbed Markov chains, Linear Algebra Appl. 417 (2006) 108–123], and the variance of the times to mixing, [J.J. Hunter, Variances of first passage times in a Markov chain with applications to mixing times, Linear Algebra Appl. 429 (2008) 1135–1162]. Some new results for the distribution of recurrence and first passage times in three-state Markov chain are also presented.

## 1. Introduction

Let $P = [p_{ij}]$ be the transition matrix of a finite irreducible, discrete time Markov chain $\{X_n\}$, $(n \geq 0)$, with state space $\mathcal{S} = \{1, 2, \ldots, m\}$.

Let $\{\pi_j\}$, $(1 \leq j \leq m)$, be the stationary distribution of the chain and $\boldsymbol{\pi}^T = (\pi_1, \pi_2, \ldots, \pi_m)$ its stationary probability vector.

For all regular (finite, aperiodic, irreducible) Markov chains, for all $j \in \mathcal{S}$, $\lim_{n \to \infty} P[X_n = j] = \pi_j$.

For all irreducible chains (including periodic chains), if for some $k \geq 0$, $P[X_k = j] = \pi_j$ for all $j \in S$, then $P[X_n = j] = \pi_j$ for all $n \geq k$ and all $j \in \mathcal{S}$.

Once the Markov chain "achieves stationarity", at say step $n$, the distribution of $X_n$ is assumed to be the stationary distribution, i.e. $P[X_n = j] = \pi_j$ for each $j \in \mathcal{S}$. If that is the case, then it easy to show that, for all $k \geq n$, $P[X_k = j] = \pi_j$ for each $j \in \mathcal{S}$.

Let $T_{ij}$ be the *"first passage time"* random variable from state $i$ to state $j$, i.e. $T_{ij} = min\{n \geq 1$ such that $X_n = j$ given that $X_0 = i\}$. Let $T_{ij}^+$ be the *"first hitting time"* random variable from state $i$ to state $j$, i.e. $T_{ij}^+ = min\{n \geq 0$ such that $X_n = j$ given that $X_0 = i\}$.

This distinction between first passage times and hitting times is only of interest when $i = j$, in which case $T_{ii}^+ = 0$ while $T_{ii} \geq 1$. For $i \neq j$, $T_{ij}^+ = T_{ij}$.

The *"mixing time"* or *"time to stationarity"* in a finite irreducible discrete time Markov chain, starting in state $i$, can be regarded as the number of trials, (the time), for the chain to reach a state sampled from the stationary distribution of the Markov chain. To be more specific:

**Definition 1**:
Let $\{X_n, n \geq 0\}$ be a Markov chain with state space $\mathcal{S} = \{1, 2, \ldots, m\}$. The random variable $M$ is said to be a *"mixing variable"* if $P[M = j] = \pi_j$ for all $j \in \mathcal{S}$, where $\{\pi_j\}$ is the stationary distribution of the chain If, under such a sampling, $M = j$, state $j$ is said to be *"the mixing state"*.

Thus the mixing state is sampled from the stationary distribution of the Markov chain.

**Definition 2:**
We say that the Markov chain $\{X_i\}$ *"achieves mixing at time $T = k$"* when $X_k = M$, the mixing variable, for the smallest such $k$.

When the concept of "mixing" was introduced in [5] the Markov chain was required to make at least a single step so that in Definition 2, $k \geq 1$, implying that mixing was achieved following a *"first passage"* from the initial state $i$ to the mixing state $j$, (or *"first return"* to state $i$ if $i = j$). However, in [3], it was found useful to permit the mixing process to terminate initially (when in Definition 2, $k = 0$) if the mixing state is the same as the initial state $i$, so that mixing occurs at the *"hitting"* time of the mixing state.

We distinguish between these two cases.

**Definition 3:**
Let $\{X_n, n \geq 0\}$ be a Markov chain with stationary distribution $\{\pi_j\}$ and mixing state $M$.
The random variable $T_i^{(0)} (\geq 0)$ is the number of trials $n$ ($n \geq 0$), given the starting (or initial) state $X_0 = i$, for the Markov chain to make a *"first hitting"* of the mixing state $M$. The random variable $T_i^{(1)} (\geq 1)$ is the number of trials $n$ ($n \geq 1$), given $X_0 = i$, for the Markov chain to make a *"first passage"* to the mixing state $M$.

While both $T_i^{(0)}$ and $T_i^{(1)}$ are *"mixing times"* of the Markov chain, starting in state $i$, we can distinguish between the two random variables by calling $T_i^{(0)}$ the *"random hitting time"* starting in state $i$ and $T_i^{(1)}$ the "*random first passage time*" starting in state $i$.

These random variables have also been used in the past as possible "mixing" variables (see [1], [9]).

Under finite state space and irreducibility conditions, the first passage times $T_{ij}$ are proper variables with finite expectations, (Theorem 7.3.1, [7]). Let $m_{ij}$ be the mean first passage time from state $i$ to state $j$, i.e. $m_{ij} = E[T_{ij} \mid X_0 = i]$ for all $i, j \in \mathcal{S}$.

Under the same conditions, the mixing times are also finite (a.s) with finite expectations. Expressions for the expected time to mixing, starting in state $i$ were derived in [5] where it was shown that $\eta_i = E\left[T_i^{(1)}\right] = \sum_{j=1}^{m} m_{ij} \pi_j$, while in [3] it was shown that $\tau_i = E\left[T_i^{(0)}\right]$ $= \sum_{j=1, j \neq i}^{m} m_{ij} \pi_j$. Thus these expectations depends on the stationary distribution of the Markov chain and the mean first passage times from state $i$ to the other states in the state space. Of considerable significance is that it was shown that these expectations are constant and neither depends on the starting state $i$, so that $\eta_i = \eta$ and further that $\tau_i = \tau = \eta - 1$. In paper [5] the main properties of $\eta$ were explored, including calculation techniques and uniform lower bounds on this expectation for all finite state Markov chains. These were extended in paper [3] to the expectation $\tau$.

In [4], expressions for the variance of the mixing times were obtained but these expressions, in general, depend on the starting state $i$.

In presenting the aforementioned results at a recent conference, the question was raised regarding the feasibility of deriving the distribution times of the mixing times $T_i^{(0)}$ and $T_i^{(1)}$. This paper provides techniques for such derivations and, further, re-establishes the expectation results above, but with different proofs. The general theory is illustrated through a study of the special cases of 2-state and 3-state Markov chains. Subsidiary to the main thrust of the paper are some new general expressions for the distributions of the first passage time and recurrence time distributions for states in a general three-state Markov chain.

## 2. Distribution Results

Let $f_{i,n} = P\{T_i^{(0)} = n \mid X_0 = i\}$ and $g_{i,n} = P\{T_i^{(1)} = n \mid X_0 = i\}$ so that $\{f_{i,n}\}$ and $\{g_{i,n}\}$ are the probability distributions of the mixing time random variables, respectively $T_i^{(0)}$ and $T_i^{(1)}$, given that the Markov chain starts in state $i$.

The $n$-step first passage time probabilities of the Markov chain $\{X_n\}$ are given as $f_{ij}^{(n)} = P\{X_n = j, X_k \neq j \text{ for } k = 1, 2, \dots n-1 \mid X_0 = i\}$, $(i, j) \in \mathcal{S} = \{1, 2, \dots, m\}$.

**Theorem 2.1**: (Distribution of the mixing times $T_i^{(0)}$ and $T_i^{(1)}$)

$$f_{i,n} = \begin{cases} \sum_{j \neq i} f_{ij}^{(n)} \pi_j, & n \geq 1, \\ \pi_i & n = 0; \end{cases} \quad (2.1)$$

and

$$g_{i,n} = \begin{cases} \sum_j f_{ij}^{(n)} \pi_j, & n \geq 1, \\ 0 & n = 0. \end{cases} \quad (2.2)$$

**Proof**: Let us assume that $X_0 = i$, so that it is given that the starting state is $i$.
First observe that $f_{i,0} = P\{T_i^{(0)} = 0\} = \sum_j P\{T_i^{(0)} = 0 | M = j\} P\{M = j\}$.

But, $\{T_i^{(0)} = 0\} \equiv \{M = i\}$ so that

$f_{i,0} = P\{T_i^{(0)} = 0 | M = i\} \pi_i + \sum_{j \neq i} P\{T_i^{(0)} = 0 | M = j\} \pi_j = \pi_i$.

In general, for $n \geq 1$, $f_{i,n} = P\{T_i^{(0)} = n\} = \sum_j P\{T_i^{(0)} = n | M = j\} P\{M = j\} = \sum_{j \neq i} f_{ij}^{(n)} \pi_j$,

since if $j = i$ mixing has occurred at the initial trial. If the mixing state is $j$ and the starting state is $i$, where $j \neq i$, then mixing can only occur for the first-time in $n$ steps if there is a first passage from state $i$ to state $j$ in $n$ steps, leading to Eqn. (2.1).

For $T_i^{(1)}$ the mixing random variable is always $\geq 1$ so that $g_{i,0} = 0$. As before, assuming that $X_0 = i$, if the mixing state is $j$, then mixing can only occur for the first-time in $n$ steps if there is a first passage from state $i$ to state $j$ in $n$ steps, (or a first return when $i = j$). i.e. $g_{i,n} = P\{T_i^{(1)} = n\} = \sum_j P\{T_i^{(1)} = n | M = j\} P\{M = j\} = \sum_j f_{ij}^{(n)} \pi_j$ leading to Eqn. (2.2).

While it is possible to use Eqns. (2.1) and (2.2) to evaluate the distributions of $T_i^{(0)}$ and $T_i^{(1)}$, these expressions require the determination of the first passage time distribution times, when typically we only have the structure of the transition matrix, $P$, and the transition probabilities, $p_{ij}$. Techniques for finding these first passage time probabilities are given in Section 5.1 of [6] and 6.2 of [7]. We do not go into these derivations in this paper, but refer the reader to the given references regarding such techniques.

Equations (2.1) and (2.2) are amenable to generating function techniques.

Let us define the probability generating functions $f_i(s) = \sum_{n=0}^{\infty} f_{i,n} s^n$ and $g_i(s) = \sum_{n=0}^{\infty} g_{i,n} s^n$.

Let $F_{ij}(s) = \sum_{n=0}^{\infty} f_{ij}^{(n)} s^n$ be the probability generating function of the first passage time random variable $T_{ij}$.

**Theorem 2.2**: (Generating functions for mixing time distributions in terms of $F_{ij}(s)$)
For $i = 1, \ldots, m$

$$f_i(s) = \pi_i + \sum_{j \neq i} F_{ij}(s) \pi_j, \quad (2.3)$$

and

$$g_i(s) = \sum_{j=1}^{m} \pi_j F_{ij}(s). \quad (2.4)$$

**Proof:** Firstly, from Eqn. (2.1),

$f_i(s) = \sum_{n=0}^{\infty} f_{i,n} s^n = \pi_i + \sum_{n=1}^{\infty} \sum_{j \neq i} f_{ij}^{(n)} s^n \pi_j = \pi_i + \sum_j \sum_{n=1}^{\infty} f_{ij}^{(n)} s^n \pi_j - \sum_{n=1}^{\infty} f_{ii}^{(n)} s^n \pi_i$

$$= \pi_i + \sum_j F_{ij}(s)\pi_j - F_{ii}(s)\pi_i, \tag{2.5}$$

leading to Eqn. (2.3).
Secondly, from Eqn. (2.2),
$$g_i(s) = \sum_{n=1}^{\infty} g_{i,n} s^n = \sum_{n=1}^{\infty} \sum_j f_{ij}^{(n)} s^n \pi_j = \sum_j \sum_{n=1}^{\infty} f_{ij}^{(n)} s^n \pi_j = \sum_j F_{ij}(s)\pi_j$$
giving Eqn. (2.4).

Let us define the n × 1 column vectors $\mathbf{f}(s) = \begin{bmatrix} f_1(s) \\ f_2(s) \\ ... \\ f_m(s) \end{bmatrix}$ and $\mathbf{g}(s) = \begin{bmatrix} g_1(s) \\ g_2(s) \\ ... \\ g_m(s) \end{bmatrix}$.

Further, define the matrix generating function $\mathbf{F}(s) = \left[ F_{ij}(s) \right]$.

**Theorem 2.3**: (Vector generating functions for mixing time distributions in terms of $\mathbf{F}(s)$)
For $|s| < 1$
$$\mathbf{f}(s) = \left[ \mathbf{F}(s) + I - \mathbf{F}_d(s) \right]\pi, \tag{2.6}$$
$$\mathbf{g}(s) = \mathbf{F}(s)\pi. \tag{2.7}$$

**Proof**: Expressing Eqn. (2.3) in vector form yields
$\mathbf{f}(s) = (I - \mathbf{F}_d(s))\pi + \mathbf{F}(s)\pi = \left[ I + \mathbf{F}(s) - \mathbf{F}_d(s) \right]\pi$, leading to Eqn. (2.6).
Similarly, Eqn. (2.7) follows directly from Eqn. (2.4).

In order to implement the results of Theorem 2.3 we need to be able to develop expressions for $\mathbf{F}(s)$ from the properties of the Markov chain. The following results provide a connection, utilizing results for the *n*-step transition probabilities $p_{ij}^{(n)} = P\{X_n = j | X_0 = i\}$.

**Theorem 2.4**: (Matrix generating function of the *n*-step transition probabilities and the *n*-step first passage time probabilities)
Let $\mathbf{P}(s) = \left[ \sum_{n=0}^{\infty} p_{ij}^{(n)} s^n \right] = \left[ P_{ij}(s) \right] = \sum_{n=0}^{\infty} P^n s^n$, then for $|s| < 1$,
$$\mathbf{P}(s) = \left[ I - Ps \right]^{-1}, \tag{2.8}$$
and
$$\mathbf{F}(s) = \left[ \mathbf{P}(s) - I \right]\left[ \mathbf{P}_d(s) \right]^{-1}, \tag{2.9}$$

where $\mathbf{P}_d(s)$ is the matrix of diagonal elements of $\mathbf{P}(s)$.

**Proof**: Eqn (2.8) is given in Theorem 6.1.9, [7] and Eqn. (2.9) is given in Theorem 6.2.6, [7].

**Theorem 2.5**: (Vector generating functions for mixing time distributions in terms of $\mathbf{P}(s)$)
For $|s| < 1$
$$\mathbf{f}(s) = \mathbf{P}(s)\left[ \mathbf{P}_d(s) \right]^{-1}\pi, \tag{2.10}$$
$$\mathbf{g}(s) = \left[ \mathbf{P}(s) - I \right]\left[ \mathbf{P}_d(s) \right]^{-1}\pi. \tag{2.11}$$

**Proof**: From Eqn. (2.9), $\mathbf{F}(s)\mathbf{P}_d(s) = \mathbf{P}(s) - I$, so that taking diagonal elements yields $\mathbf{F}_d(s)\mathbf{P}_d(s) = \mathbf{P}_d(s) - I$ implying $\mathbf{F}_d(s) = I - [\mathbf{P}_d(s)]^{-1}$ and $I - \mathbf{F}_d(s) = [\mathbf{P}_d(s)]^{-1}$. Eqn. (2.10) follows from Eqn. (2.6), while Eqn. (2.11) follows directly from Eqns. (2.7) and (2.9).

From the above results one notes that elemental expressions for the generating functions $f_i(s)$ and $g_i(s)$ can be given using equations (2.10) and (2.11), respectively.

**Theorem 2.6**: (Generating functions for mixing time distributions in terms of $P_{ij}(s)$)
For $|s| < 1$

$$f_i(s) = \sum_{j=1}^{m} \pi_j \left[ \frac{P_{ij}(s)}{P_{jj}(s)} \right], \tag{2.12}$$

$$g_i(s) = f_i(s) - \frac{\pi_i}{P_{ii}(s)}. \tag{2.13}$$

Note that the results of Theorems 2.2. and 2.6 are also linked by results connecting the generating functions $F_{ij}(s)$ and $P_{ij}(s)$. From Theorem 6.2.5 of [7], for all $i, j \in \{1, 2, \ldots, m\}$,

$$F_{ij}(s) = \frac{P_{ij}(s) - \delta_{ij}}{P_{jj}(s)} \text{ and } \frac{1}{P_{ii}(s)} = 1 - F_{ii}(s). \tag{2.14}$$

Theorem 6.1.10 of [7] describes the expansion of $\mathbf{P}(s) = [I - Ps]^{-1}$ when the transition matrix $P$ has distinct eigenvalues $\lambda_1 = 1, \lambda_2, \ldots, \lambda_m$. Let $\mathbf{x}_1, \mathbf{x}_2, \ldots, \mathbf{x}_m$ and $\mathbf{y}_1^T, \mathbf{y}_2^T, \ldots, \mathbf{y}_m^T$ be the corresponding right and left eigenvectors chosen so that $\mathbf{y}_i^T \mathbf{x}_i = 1$, $(i = 1, 2, \ldots, m)$. Then for $|s| < 1$, $\mathbf{P}(s) = [I - Ps]^{-1} = \sum_{k=1}^{m} \sum_{n=1}^{\infty} A_k \lambda_k^n s^n$ where for nonzero $\lambda_k$, $A_k = \mathbf{x}_k \mathbf{y}_k^T$.

Observe that $[I - Ps]^{-1} = \sum_{k=1}^{m} \frac{A_k}{1 - \lambda_k s}$ and that $P^{(n)} = [p_{ij}^{(n)}] = \sum_{k=1}^{m} A_k \lambda_k^n$.

If we define $\mathbf{x}_k^T = (x_{1k}, x_{2k}, \ldots, x_{mk})$ and $\mathbf{y}_k^T = (y_{k1}, y_{k2}, \ldots, y_{km})$ then $A_k = [x_{ik} y_{kj}]$ so that

$$P_{ij}(s) = \sum_{k=1}^{m} \frac{x_{ik} y_{kj}}{1 - \lambda_k s}.$$

Further, $\det(I - Ps) = (1 - \lambda_1 s)(1 - \lambda_2 s) \ldots (1 - \lambda_m s)$ where $\lambda_1 = 1, \lambda_2, \ldots, \lambda_m$ are the eigenvalues of $P$ and that $A_k$ can be found directly as $A_k = \frac{\lambda_k^{m-1}}{\prod_{j \neq k}(\lambda_k - \lambda_j)} \text{adj}(I - \lambda_k^{-1} P) = \mathbf{x}_k \mathbf{y}_k^T$.

Note that the characteristic polynomial of P is $c(\lambda) = \det(P - \lambda I) = (-1)^m \lambda^m \det(I - \lambda^{-1} P)$.

Further $\mathbf{P}(s) = [I - Ps]^{-1} = [P_{ij}(s)] = \frac{1}{\det(I - Ps)} \text{adj}(I - Ps) = \frac{1}{\det(I - Ps)} [a_{ij}(s)]$, (2.15)

where $a_{ij}(s)$ is the $(j,i)$-th cofactor of $I - sP$.

Now $\mathbf{P}_d(s) = \frac{1}{\det(I - sP)} [\delta_{ij} a_{jj}(s)]$ so that $[\mathbf{P}_d(s)]^{-1} = \det(I - sP) \left[ \frac{\delta_{ij}}{a_{jj}(s)} \right]$.

Consequently $\mathbf{P}(s)\left[\mathbf{P}_d(s)\right]^{-1} = \left[\dfrac{P_{ij}(s)}{P_{jj}(s)}\right] = \left[\dfrac{a_{ij}(s)}{a_{jj}(s)}\right]$ and thus when one wishes to evaluate the ratio $P_{ij}(s)/P_{jj}(s)$ one does not need to compute the determinant, only the elements of the adjoint.

The implementation of these results is best illustrated in some examples. See Sections 3 and 4.

While general expressions for the distributions of the mixing times are difficult to obtain, it is relatively easy to extract moments of the mixing times, using results of the moments of the first passage times and the relationships given by Eqns. (2.3) and (2.4).

**Theorem 2.7:** (Expected times to mixing)
If the mean first passage time from state $i$ to state $j$ is $m_{ij}$ then
$$E\left[T_i^{(0)}\right] = \sum_{j \neq i} \pi_j m_{ij}, \tag{2.16}$$
$$E\left[T_i^{(1)}\right] = \sum_{j} \pi_j m_{ij}. \tag{2.17}$$

**Proof**: Since $E\left[T_i^{(0)}\right] = \lim_{s \uparrow 1} \dfrac{df_i(s)}{ds}$ and since from Eqn. (2.3), $\dfrac{df_i(s)}{ds} = \sum_{j \neq i} \dfrac{dF_{ij}(s)}{ds} \pi_j$, taking the limit as $s \uparrow 1$, and noting that $\lim_{s \uparrow 1} \dfrac{F_{ij}(s)}{ds} = m_{ij}$ yields Eqn. (2.16). Similarly Eqn. (2.17) follows from Eqn. (2.4).

The fact that these means are invariant under changing the initial starting state $i$ is a curious phenomena. The derivation of the result that $\sum_j \pi_j m_{ij} = \eta$ are independent of $i$ is discussed in [5], and the result that $\sum_{j \neq i} \pi_j m_{ij} = \tau = \eta - 1$ is discussed in [3]. The linking of the two expectations follows from the observation that $m_{ii} \pi_i = 1$. We do not repeat the derivation of these results but note that various expressions for $\eta$ and $\tau$ can be given, typically involving the trace of generalized inverses of $I - P$. In particular, $\eta = tr(Z)$, where $Z = [I - P + \Pi]^{-1}$, (with $\Pi = e\pi^T$). $Z$ is Kemeny and Snell's fundamental matrix, ([8]). The constant $\eta$ is also known as Kemeny's constant (see [2, Chapter 11], [8, Corollary 4.3.6]).

In [5] it is shown that for irreducible periodic, period $m$, Markov chains, $\eta = (m+1)/2$; for an $m$-state Markov chain consisting of independent trials, $\eta = m$; while for any irreducible $m$-state Markov chain, $\eta \geq (m+1)/2$.

## 3. Special case –Two-state Markov chains

Let $P = \begin{bmatrix} p_{11} & p_{12} \\ p_{21} & p_{22} \end{bmatrix} = \begin{bmatrix} 1-a & a \\ b & 1-b \end{bmatrix}$ (3.1)

with $0 \leq a \leq 1$, $0 \leq b \leq 1$, be the transition matrix of a two-state Markov chain with state space $S = \{1, 2\}$. Let $d = 1 - p_{12} - p_{21} = 1 - a - b$.

If $-1 \leq d < 1$, the Markov chain is irreducible with a unique stationary distribution given by

$$\pi_1 = \frac{p_{21}}{p_{12} + p_{21}} = \frac{b}{1-d}, \quad \pi_2 = \frac{p_{12}}{p_{12} + p_{21}} = \frac{a}{1-d}. \tag{3.2}$$

If $-1 < d < 1$, the Markov chain is regular and this stationary distribution is in fact the limiting distribution. If $d = 1$, there is no unique stationary distribution (with both states absorbing), while if $d = -1$ the Markov chain is irreducible periodic, period 2.

In the case of independent trials, $P = \begin{bmatrix} p_{11} & p_{12} \\ p_{11} & p_{12} \end{bmatrix} = \begin{bmatrix} 1-a & a \\ 1-a & a \end{bmatrix}$, with identical rows so that $b = 1 - a$ and $d = 0$.

For this two-state Markov chain, from Example 6.1.6 [7],

$$\mathbf{P}(s) = \frac{1}{(1-s)(1-ds)} \begin{bmatrix} 1-p_{22}s & p_{12}s \\ p_{21}s & 1-p_{11}s \end{bmatrix} = \frac{1}{(1-s)(1-ds)} \begin{bmatrix} 1-(1-b)s & as \\ bs & 1-(1-a)s \end{bmatrix}, \tag{3.3}$$

while, from Exercise 6.2.2 [7],

$$\mathbf{F}(s) = \begin{bmatrix} \dfrac{s(p_{11}-ds)}{1-p_{22}s} & \dfrac{p_{12}s}{1-p_{11}s} \\ \dfrac{p_{21}s}{1-p_{22}s} & \dfrac{s(p_{22}-ds)}{1-p_{11}s} \end{bmatrix} = \begin{bmatrix} \dfrac{s(1-a-ds)}{1-(1-b)s} & \dfrac{as}{1-(1-a)s} \\ \dfrac{bs}{1-(1-b)s} & \dfrac{s(1-b-ds)}{1-(1-a)s} \end{bmatrix}. \tag{3.4}$$

We first summarise the results for the distribution of the recurrence time r.v. $T_{11}$ and the first passage time r.v. $T_{11}$.

**Theorem 3.1:** (The distributions of $T_{11}$ and $T_{12}$ for 2-state Markov chains)

$$f_{11}^{(1)} = p_{11}, \ f_{11}^{(n)} = p_{12}p_{22}^{n-2}p_{21}, \ (n \geq 2), \tag{3.5}$$

$$f_{12}^{(n)} = p_{11}^{n-1}p_{12}, \ (n \geq 1). \tag{3.6}$$

**Proof**: The proofs are well known (Theorem 5.1.8, [6]) and follow from extracting the coefficient of $s^n$ from $F_{11}(s)$ and $F_{12}(s)$, as given in Eqn.(3.4). Alternatively Eqns. (3.5) and (3.6) follow by using simple sample path arguments.

Now from Eqn. (3.3)

$$\mathbf{P}_d(s) = \frac{1}{(1-s)(1-ds)} \begin{bmatrix} 1-p_{22}s & 0 \\ 0 & 1-p_{11}s \end{bmatrix},$$

so that 
$$[\mathbf{P}_d(s)]^{-1} = (1-s)(1-ds) \begin{bmatrix} \dfrac{1}{1-p_{22}s} & 0 \\ 0 & \dfrac{1}{1-p_{11}s} \end{bmatrix}.$$

From Eqn. (2.9) and Eqn. (3.4),

$$\mathbf{f}(s) = \begin{bmatrix} f_1(s) \\ f_2(s) \end{bmatrix} = \begin{bmatrix} 1-p_{22}s & p_{12}s \\ p_{21}s & 1-p_{11}s \end{bmatrix} \cdot \begin{bmatrix} \dfrac{1}{1-p_{22}s} & 0 \\ 0 & \dfrac{1}{1-p_{11}s} \end{bmatrix} \begin{bmatrix} \pi_1 \\ \pi_2 \end{bmatrix} = \begin{bmatrix} \pi_1 + \pi_2 \dfrac{p_{12}s}{1-p_{11}s} \\ \pi_2 + \pi_1 \dfrac{p_{21}s}{1-p_{22}s} \end{bmatrix}. \qquad (3.7)$$

From Eqn. (2.4) and Eqn. (3.3),

$$\mathbf{g}(s) = \begin{bmatrix} g_1(s) \\ g_2(s) \end{bmatrix} = \begin{bmatrix} \dfrac{s(p_{11}-ds)}{1-p_{22}s} & \dfrac{p_{12}s}{1-p_{11}s} \\ \dfrac{p_{21}s}{1-p_{22}s} & \dfrac{s(p_{22}-ds)}{1-p_{11}s} \end{bmatrix} \begin{bmatrix} \pi_1 \\ \pi_2 \end{bmatrix} = \begin{bmatrix} \pi_1 \dfrac{s(p_{11}-ds)}{1-p_{22}s} + \pi_2 \dfrac{p_{12}s}{1-p_{11}s} \\ \pi_1 \dfrac{p_{21}s}{1-p_{22}s} + \pi_2 \dfrac{s(p_{22}-ds)}{1-p_{11}s} \end{bmatrix}. \qquad (3.8)$$

Note that for all cases where $-1 \le d < 1$, $f_1(1) = f_2(1) = \pi_1 + \pi_2 = 1$ and $g_1(1) = g_2(1) = \pi_1 + \pi_2 = 1$, reconfirming that $T_i^{(0)}$ and $T_i^{(1)}$ are both proper random variables.

**Theorem 3.2:** (The distributions of $T_1^{(0)}$ and $T_1^{(1)}$ for 2-state Markov chains)
For the two-state Markov chain with transition matrix given by Eqn. (3.1), if $-1 \le d < 1$, the distribution of the mixing time random variable $T_1^{(0)}$ is given by

$$f_{10} = \pi_1, f_{1n} = \pi_2 p_{12} p_{11}^{n-1}, (n \ge 1), \qquad (3.9)$$

where the stationary distribution $\{\pi_j\}$ is given by Eqn. (3.2).
If $-1 \le d < 1$, the distribution of the mixing time random variable $T_1^{(1)}$ is given by

$$g_{10} = 0, \; g_{11} = \pi_1 p_{11} + \pi_2 p_{12}, \; g_{1n} = \pi_1(p_{11}p_{22} - d)p_{22}^{n-2} + \pi_2 p_{12} p_{11}^{n-1}, (n \ge 2). \qquad (3.10)$$

**Proof**: Expanding the power series for $f_1(s)$, given in Eqn. (3.5), we obtain, for starting in state 1, (with symmetrical results for starting in state 2),

$$f_1(s) = \pi_1 + \pi_2 \frac{p_{12}s}{1-p_{11}s} = \pi_1 + \pi_2 \sum_{k=0}^{\infty} p_{12} p_{11}^k s^{k+1}$$ leading to expression given by Eqn. (3.9).

Alternatively, Eqn. (3.9) follows directly from Eqn. (2.1) and Eqns. (3.5) and (3.6).
Similarly, expanding the power series for $g_1(s)$, given in Eqn. (3.8),

$$g_1(s) = \pi_1 \frac{s(p_{11}-ds)}{1-p_{22}s} + \pi_2 \frac{p_{12}s}{1-p_{11}s}$$

$$= \pi_1 p_{11} \sum_{k=0}^{\infty} p_{22}^k s^{k+1} - \pi_1 d \sum_{k=0}^{\infty} p_{22}^k s^{k+2} + \pi_2 p_{12} \sum_{k=0}^{\infty} p_{11}^k s^{k+1}$$

leading to the expressions given in Eqn. (3.10).
Alternatively, Eqn. (3.10) follows directly from Eqn. (2.2) and Eqns. (3.5) and (3.6).

Theorem 3.2 establishes that $T_1^{(0)}$ is a modified geometric random variable and that $T_1^{(1)}$ is the mixture of two geometric random variables. Similar results hold for $T_2^{(0)}$ and $T_2^{(1)}$ (by interchanging 1 and 2).

The expected times to mixing can be obtained from the generating functions $f_1(s)$ and $g_1(s)$.

**Theorem 3.3:** (Mean mixing times for 2-state Markov chains)
For the two-state Markov chain with transition matrix given by Eqn. (3.1), if $-1 \leq d < 1$,

$$E\left[T_1^{(0)}\right] = E\left[T_2^{(0)}\right] = \frac{1}{1-d} = \tau, \qquad (3.11)$$

and
$$E\left[T_1^{(1)}\right] = E\left[T_2^{(1)}\right] = 1 + \frac{1}{1-d} = \eta. \qquad (3.12)$$

**Proof:** From Eqn. (3.7), since $f_1(s) = \pi_1 + \pi_2 \frac{p_{12} s}{1 - p_{11} s}$,

$$\frac{df_1(s)}{ds} = \pi_2 p_{12} \frac{d}{ds} s[1 - p_{11} s]^{-1} = \pi_2 p_{12} \left[ (1 - p_{11} s)^{-1} - s(1 - p_{11} s)^{-2}(-p_{11}) \right].$$

Thus $E\left[T_1^{(0)}\right] = \lim_{s \uparrow 1} \frac{df_1(s)}{ds} = \pi_2 p_{12} \left[ (1 - p_{11})^{-1} + p_{11}(1 - p_{11})^{-2} \right]$

$$= \pi_2 p_{12} \left[ \frac{1}{p_{12}} + \frac{p_{11}}{p_{12}^2} \right] = \frac{\pi_2 p_{12}}{p_{12}^2}(p_{12} + p_{11}) = \frac{\pi_2}{p_{12}} = \frac{1}{p_{12} + p_{21}} = \frac{1}{1-d}.$$

By the symmetry of the above result, interchanging the indices 1 and 2, leads to identical expressions, as given by Eqn. (3.11).

Similarly, from Eqn. (3.8), since $g_1(s) = \pi_1 \frac{s(p_{11} - ds)}{1 - p_{22} s} + \pi_2 \frac{p_{12} s}{1 - p_{11} s}$,

$$\frac{dg_1(s)}{ds} = \pi_1 \frac{d}{ds}\left[(p_{11}s - ds^2)(1 - p_{22}s)^{-1}\right] + \pi_2 \frac{d}{ds}\left[p_{12}s(1 - p_{11}s)^{-1}\right]$$

$$= \pi_1 \left[(p_{11} - 2ds)(1 - p_{22}s)^{-1} + p_{22}(p_{11}s - ds^2)(1 - p_{22}s)^{-2}\right]$$

$$+ \pi_2 \left[p_{12}(1 - p_{11}s)^{-1} - p_{11}p_{12}s(1 - p_{11}s)^{-2}\right].$$

Now $E\left[T_1^{(1)}\right] = \lim_{s \uparrow 1} \frac{dg_1(s)}{ds} = \pi_1 \left[\frac{p_{11} - 2d}{p_{21}} + \frac{p_{22}(p_{11} - d)}{p_{21}^2}\right] + \pi_2 \left[1 + \frac{p_{11}}{p_{12}}\right],$

$$= \frac{1}{p_{12} + p_{21}}(p_{11} - 2d + p_{22} + 1), \text{ since } p_{11} - d = p_{21},$$

$$= \frac{1}{p_{12} + p_{21}}(p_{21} - d + p_{22} + 1) = \frac{1}{1-d}(2 - d) = 1 + \frac{1}{1-d}.$$

In the above proof we have established expressions for the expected times to mixing for each starting state, without resorting to the complicated arguments that were used to derive these results in a general setting in [5] (for the case of $\eta$) and in [3] (for the case of $\tau$).

For all two-state irreducible Markov chains, $\tau \geq 0.5$, ([3]), and $\eta \geq 1.5$, ([5]), with arbitrarily large values of $\tau$ and $\eta$ occurring as $d \to 1$, (when both $a \to 0$ and $b \to 0$). This occurs when the chain is approaching the situation of being close to reducible, with both states absorbing.

**Periodic Markov chains**
Note that when $d = -1$, $p_{12} = p_{21} = 1$, so that the Markov chain is periodic, period 2, with $\pi_1 = \pi_2 = 1/2$. Under these conditions the results of Theorem 3.2 yield $f_{10} = \pi_1$, $f_{11} = \pi_2$, and $f_{1n} = 0$, $(n \geq 2)$, and $g_{10} = 0$, $g_{11} = \pi_2$, $g_{12} = \pi_1$, and $g_{1n} = 0$, $(n \geq 3)$.

This is consistent with the following observations. Suppose that the Markov chain starts in state 1 with $X_0 = 1$. If the mixing state $M$ is 1 (with probability 1/2) then the random hitting time $T_1^{(0)} = 0$, so that mixing occurs at that trial while the random first passage time $T_1^{(1)} = 2$ since 2 further steps $1 \to 2, 2 \to 1$ are required. If the mixing state $M$ is 2 (with probability 1/2) then $T_1^{(0)} = T_1^{(1)} = 1$ since the mixing state occurs after 1 further step as $1 \to 2$.

The minimum value of the expected mixing times are $\tau = 0.5, \eta = 1.5$ which occur when $d = -1$, i.e. in this periodic, period 2 case.

**Independent trials**
In the case of independent trials with two outcomes (states 1 and 2), $b = 1 - a$, $d = 0$, and $\pi_1 = 1 - a, \pi_2 = a$.

From the results of Theorem 3.2, the distribution of the mixing time random variable $T_1^{(0)}$ is given by $f_{10} = 1 - a, f_{1n} = a^2(1-a)^{n-1}, (n \geq 1)$, and the distribution of the mixing time random variable $T_1^{(1)}$ is given by $g_{10} = 0, g_{11} = (1-a)^2 + a^2, g_{1n} = (1-a)^2 a^{n-1} + a^2(1-a)^{n-1}, (n \geq 2)$.

In independent trials, the mixing time $T_1^{(1)}$ is effectively the time for a nominated state (1 or 2) to occur under repeated identical conditions, so that we have the sum of two weighted geometric random variables with parameters $a$ and $1 - a$ with weights $(\pi_1 =)1 - a$ or $(\pi_2 =) a$ depending on whether we are waiting for state 1 or 2 to occur. For $T_1^{(0)}$, we either have an occurrence initially (with probability $1 - a$) or we wait (with probability $a$) for a geometric random time for the other state to occur.

In the case of independent trials, since $d = 0$, the expected times to mixing are $\tau = 1$ and $\eta = 2$.

## 4. Special case – Three-state Markov chains

Let $P = \begin{bmatrix} p_{11} & p_{12} & p_{13} \\ p_{21} & p_{22} & p_{23} \\ p_{31} & p_{32} & p_{33} \end{bmatrix}$ (4.1)

be the transition matrix of a three-state Markov chain with state space $S = \{1, 2, 3\}$.

Let $\Delta_1 \equiv p_{23}p_{31} + p_{21}p_{32} + p_{21}p_{31}, \Delta_2 \equiv p_{31}p_{12} + p_{32}p_{13} + p_{32}p_{12}, \Delta_3 \equiv p_{12}p_{23} + p_{13}p_{21} + p_{13}p_{23}$, and $\Delta \equiv \Delta_1 + \Delta_2 + \Delta_3$.

The Markov chain, with the above transition matrix, is irreducible (and hence a stationary distribution exists) if and only if $\Delta_1 > 0, \Delta_2 > 0, \Delta_3 > 0$. Under these conditions, it is easily shown that the stationary probability vector is

$$(\pi_1, \pi_2, \pi_3) = \frac{1}{\Delta}(\Delta_1, \Delta_2, \Delta_3).$$ (4.2)

Special cases of this Markov chain were considered in [7] but in no instance was a general form of **P**(s) or **F**(s) derived. We explore this now in the context of the results of this paper.

Note from Eqn.(2.8) that

$$\mathbf{P}(s) = \left[P_{ij}(s)\right] = \frac{1}{\det(I-Ps)}\left[a_{ij}(s)\right] = \frac{1}{\det(I-Ps)}\begin{bmatrix} a_{11}(s) & a_{12}(s) & a_{13}(s) \\ a_{21}(s) & a_{22}(s) & a_{23}(s) \\ a_{31}(s) & a_{32}(s) & a_{33}(s) \end{bmatrix}, \quad (4.3)$$

where $a_{ij}(s)$ is the (j,i)-th cofactor of $I - sP$.

It is easily verified that

$$\det(I-Ps) = 1 - s(p_{11}+p_{22}+p_{33}) + s^2(p_{11}p_{22}+p_{22}p_{33}+p_{33}p_{11}-p_{12}p_{21}-p_{23}p_{32}-p_{13}p_{31})$$
$$- s^3(p_{12}p_{23}p_{31}+p_{13}p_{32}p_{21}+p_{11}p_{22}p_{33}-p_{12}p_{21}p_{33}-p_{13}p_{31}p_{22}-p_{11}p_{23}p_{32}). \quad (4.4)$$

Further

$$\det(I-Ps) = (1-s)(1-\lambda_2 s)(1-\lambda_3 s) = (1-s)\{1-(\lambda_2+\lambda_3)s + \lambda_2\lambda_3 s^2\}$$
$$= (1-s)(1+As+Bs^2) = 1+(A-1)s+(B-A)s^2 - Bs^3, \quad (4.5)$$

so that $A = 1-(p_{11}+p_{22}+p_{33})$ and there are two equivalents forms of $B$, viz.

$$B = p_{12}p_{23}p_{31}+p_{13}p_{32}p_{21}+p_{11}p_{22}p_{33}-p_{12}p_{21}p_{33}-p_{13}p_{31}p_{22}-p_{11}p_{23}p_{32}$$
$$= 1-(p_{11}+p_{22}+p_{33})+p_{11}p_{22}+p_{22}p_{33}+p_{33}p_{11}-p_{12}p_{21}-p_{23}p_{32}-p_{13}p_{31}$$

Note that $1-A = p_{11}+p_{22}+p_{33}$ and $A-B = p_{12}p_{21}+p_{23}p_{32}+p_{13}p_{31}-p_{11}p_{22}-p_{22}p_{33}-p_{33}p_{11}$.

For the cofactor terms,

$$a_{ii}(s) = 1 - s(p_{aa}+p_{bb}) + s^2(p_{aa}p_{bb}-p_{ab}p_{ba}) \text{ where } i,a,b \in \{1,2,3\} \text{ with } i \neq a \neq b. \quad (4.6)$$

Note $a_{ii}(s) = (1-\lambda_{ia}s)(1-\lambda_{ib}s)$ where $\lambda_{ia}$ and $\lambda_{ib}$ are roots of the quadratic, i.e.

$$\lambda_{ia} = \frac{p_{aa}+p_{bb}+\delta_i}{2}, \lambda_{ib} = \frac{p_{aa}+p_{bb}-\delta_i}{2}, \text{ with } \delta_i = \sqrt{(p_{aa}-p_{bb})^2 + 4p_{ab}p_{ba}} \ (\geq 0).$$

Also $a_{ij}(s) = sp_{ij} + s^2(p_{ik}p_{kj}-p_{ij}p_{kk})$ where $i,j,k \in \{1,2,3\}$ with $i \neq j \neq k$. (4.7)

From Eqn. (2.14),

$$\mathbf{F}(s) = \left[F_{ij}(s)\right] = \begin{bmatrix} 1-\dfrac{1}{P_{11}(s)} & \dfrac{P_{12}(s)}{P_{22}(s)} & \dfrac{P_{13}(s)}{P_{33}(s)} \\ \dfrac{P_{21}(s)}{P_{11}(s)} & 1-\dfrac{1}{P_{11}(s)} & \dfrac{P_{23}(s)}{P_{33}(s)} \\ \dfrac{P_{31}(s)}{P_{11}(s)} & \dfrac{P_{32}(s)}{P_{22}(s)} & 1-\dfrac{1}{P_{33}(s)} \end{bmatrix},$$

where $F_{11}(s) = 1 - \dfrac{1}{P_{11}(s)} = 1 - \dfrac{\det(I-Ps)}{a_{11}(s)}, F_{12}(s) = \dfrac{a_{12}(s)}{a_{22}(s)}$ and $F_{13}(s) = \dfrac{a_{13}(s)}{a_{33}(s)}$. (4.8)

From the expressions above for $F_{11}(s)$ and $F_{12}(s)$, we can derive expressions for the distributions of the recurrence time distribution of $T_{11}$ and the first passage time distribution of $T_{12}$. By symmetry, expressions for the distributions of the other $T_{ij}$ can also be obtained.

**Theorem 4.1:** (The distribution of $T_{11}$ for 3-state Markov chains).
The probability distribution $\{f_{11}^{(n)}\}$ where $f_{11}^{(n)} = P\{T_{11} = n\}$ is given by
$$f_{11}^{(1)} = p_{11}, f_{11}^{(2)} = p_{12}p_{21} + p_{13}p_{31}, \quad f_{11}^{(3)} = p_{12}p_{22}p_{21} + p_{12}p_{23}p_{31} + p_{13}p_{32}p_{21} + p_{13}p_{33}p_{31}$$
and, in general, for $n \geq 3$, provided $\delta_1 > 0$,
$$f_{11}^{(n)} = \frac{c(\lambda_{12})}{\lambda_{12} - \lambda_{13}} \lambda_{12}^{n-2} - \frac{c(\lambda_{13})}{\lambda_{12} - \lambda_{13}} \lambda_{13}^{n-2}, \tag{4.9}$$
where $c(\lambda) = \det(P - \lambda I)$, the characteristic polynomial, with
$$c(\lambda_{12}) = (1 - p_{11})(1 - p_{22})(1 - p_{33}) - p_{23}p_{32}(1 - p_{11}) - (p_{12}p_{21} + p_{13}p_{31})(1 - \lambda_{12}), \tag{4.10}$$
$$c(\lambda_{13}) = (1 - p_{11})(1 - p_{22})(1 - p_{33}) - p_{23}p_{32}(1 - p_{11}) - (p_{12}p_{21} + p_{13}p_{31})(1 - \lambda_{13}), \tag{4.11}$$
and $\lambda_{12} = \dfrac{p_{22} + p_{33} + \delta_1}{2}, \lambda_{13} = \dfrac{p_{22} + p_{33} - \delta_1}{2}$ with $\delta_1 = \sqrt{(p_{22} - p_{33})^2 + 4p_{23}p_{32}}$. \hfill (4.12)

**Proof:** We use Eqn.(4.8) to determine $F_{11}(s)$, with $\det(I - Ps)$ as given by Eqns.(4.4) and (4.5). Now $a_{11}(s) = 1 - s(p_{22} + p_{33}) + s^2(p_{22}p_{33} - p_{23}p_{32}) = (1 - \lambda_{12}s)(1 - \lambda_{13}s)$,
where $\lambda_{12} + \lambda_{13} = p_{22} + p_{33}$ and $\lambda_{12}\lambda_{13} = p_{22}p_{33} - p_{23}p_{32}$, \hfill (4.13)
leading to the terms given by Eqn.(4.12).
Thus $\dfrac{1}{a_{11}(s)} = \dfrac{1}{(1 - \lambda_{12}s)(1 - \lambda_{13}s)} = \dfrac{C}{(1 - \lambda_{12}s)} + \dfrac{D}{(1 - \lambda_{13}s)} = \sum_{n=o}^{\infty}(C\lambda_{12}^n + D\lambda_{13}^n)s^n$,
where , (since $C + D = 1$, $C\lambda_{13} + D\lambda_{12} = 0$), $C = \dfrac{\lambda_{12}}{\lambda_{12} - \lambda_{13}}$ and $D = \dfrac{-\lambda_{13}}{\lambda_{12} - \lambda_{13}}$. \hfill (4.14)

Now, from Eqns.(4.8) and (4.5),
$$F_{11}(s) = \sum_{n=0}^{\infty} f_{11}^{(n)} s^n = 1 + \{-1 + (1 - A)s + (A - B)s^2 + Bs^3\}\{\sum_{n=o}^{\infty}(C\lambda_{12}^n + D\lambda_{13}^n)s^n\}. \tag{4.15}$$
Using the results of Eqns. (4.14), observe that
$$C\lambda_{12} + D\lambda_{13} = \left(\frac{\lambda_{12}^2 - \lambda_{13}^2}{\lambda_{12} - \lambda_{13}}\right) = \lambda_{12} + \lambda_{13} = p_{22} + p_{33},$$
$$C\lambda_{12}^2 + D\lambda_{12}^2 = \left(\frac{\lambda_{12}^3 - \lambda_{13}^3}{\lambda_{12} - \lambda_{13}}\right) = \lambda_{12}^2 + \lambda_{13}^2 + \lambda_{12}\lambda_{13} = (\lambda_{12} + \lambda_{13})^2 - \lambda_{12}\lambda_{13} = p_{22}^2 + p_{33}^2 + p_{22}p_{33} + p_{23}p_{32},$$
$$C\lambda_{12}^3 + D\lambda_{12}^3 = \left(\frac{\lambda_{12}^4 - \lambda_{13}^4}{\lambda_{12} - \lambda_{13}}\right) = \left(\frac{(\lambda_{12}^2 - \lambda_{13}^2)(\lambda_{12}^2 + \lambda_{13}^2)}{\lambda_{12} - \lambda_{13}}\right) = (\lambda_{12} + \lambda_{13})\{(\lambda_{12} + \lambda_{13})^2 - 2\lambda_{12}\lambda_{13}\}$$
$$= (p_{22} + p_{33})(p_{22}^2 + p_{33}^2 + 2p_{23}p_{32}).$$
Equating the coefficients of $s^n$ for $n = 0, 1, 2$ and 3, and using the above results we obtain:
$$f_{11}^{(0)} = \text{coeff of } s^0 = 1 - (C + D) = 0,$$
$$f_{11}^{(1)} = \text{coeff of } s^1 = -(C\lambda_{12} + D\lambda_{13}) + (1 - A)(C + D) = -(p_{22} + p_{33}) + p_{11} + p_{22} + p_{33} = p_{11},$$

$$f_{11}^{(2)} = \text{coeff of } s^2 = -(C\lambda_{12}^2 + D\lambda_{13}^2) + (1-A)(C\lambda_{12} + D\lambda_{13})_1 + (A-B)(C+D)$$
$$= -(p_{22}^2 + p_{33}^2 + p_{22}p_{33} + p_{23}p_{32}) + (p_{11} + p_{22} + p_{33})(p_{22} + p_{33})$$
$$+ (p_{12}p_{21} + p_{23}p_{32} + p_{13}p_{31} - p_{11}p_{22} - p_{22}p_{33} - p_{33}p_{11}) = p_{12}p_{21} + p_{13}p_{31},$$
$$f_{11}^{(3)} = \text{coeff of } s^3 = -(C\lambda_{12}^3 + D\lambda_{13}^3) + (1-A)(C\lambda_{12}^2 + D\lambda_{13}^2) + (A-B)(C\lambda_{12} + D\lambda_{13}) + B(C+D)$$
$$= (p_{22} + p_{33})(p_{22}^2 + p_{33}^2 + 2p_{23}p_{32}) + (p_{11} + p_{22} + p_{33})(p_{22}^2 + p_{33}^2 + p_{22}p_{33} + p_{23}p_{32})$$
$$+ (p_{12}p_{21} + p_{23}p_{32} + p_{13}p_{31} - p_{11}p_{22} - p_{22}p_{33} - p_{33}p_{11})(p_{22} + p_{33})$$
$$+ p_{12}p_{23}p_{31} + p_{13}p_{32}p_{21} + p_{11}p_{22}p_{33} - p_{12}p_{21}p_{33} - p_{13}p_{31}p_{22} - p_{11}p_{23}p_{32}$$
$$= p_{12}p_{22}p_{21} + p_{12}p_{23}p_{31} + p_{13}p_{32}p_{21} + p_{13}p_{33}p_{31},$$

leading to the special cases when $n = 1, 2,$ and 3.

For the general case of the Theorem, when $n \geq 3$,
$$f_{11}^{(n)} = \text{coeff of } s^n = -(C\lambda_{12}^n + D\lambda_{13}^n) + (1-A)(C\lambda_{12}^{n-1} + D\lambda_{13}^{n-1}) + (A-B)(C\lambda_{12}^{n-2} + D\lambda_{13}^{n-2}) + B(C\lambda_{12}^{n-3} + D\lambda_{13}^{n-3})$$
$$= C\lambda_{12}^{n-3}\{-\lambda_{12}^3 + (1-A)\lambda_{12}^2 + (A-B)\lambda_{12} + B\} + D\lambda_{13}^{n-3}\{-\lambda_{13}^3 + (1-A)\lambda_{13}^2 + (A-B)\lambda_{13} + B\}$$
$$= C\lambda_{12}^{n-3}c(\lambda_{12}) + D\lambda_{13}^{n-3}c(\lambda_{13}),$$

which reduces to expression (4.9), using Eqns. (4.12). Observe further, from Eqn. (4.5), $c(\lambda) = -\lambda^3 \det(I - \lambda^{-1}P) = -\lambda^3 + (1-A)\lambda^2 + (A-B)\lambda + B$, the characteristic polynomial. Now $c(\lambda) = (1-\lambda)(B + A\lambda + \lambda^2)$, and substitution by $\lambda = \lambda_{12}$ and $\lambda = \lambda_{13}$, using Eqns. (4.12), yields, after simplification, the expressions given by Eqns. (4.10) and (4.11).
Note also that the expression given for $f_{11}^{(3)}$ also follows from Eqn. (4.9) when $n = 3$.

An immediate observation is that the recurrence time distribution for state 1 (and similarly for the other states) appears as a mixture of two geometric distributions, although as we see this can reduce to a single geometric distribution. In the case where $\delta_1 = 0$, the distribution can reduce to a negative binomial (see Case 3 to follow.) An extension to the above result is a recurrence relationship that can be used as an alternative computational procedure.

**Corollary 4.2:** (The distribution of $T_{11}$ for 3-state Markov chains).
The probability distribution $\{f_{11}^{(n)}\}$ where $f_{11}^{(n)} = P\{T_{11} = n\}$ is given by
$$f_{11}^{(1)} = p_{11}, f_{11}^{(2)} = p_{12}p_{21} + p_{13}p_{31}, f_{11}^{(3)} = p_{12}p_{22}p_{21} + p_{12}p_{23}p_{31} + p_{13}p_{32}p_{21} + p_{13}p_{33}p_{31},$$
and, in general, for $n \geq 4$,
$$f_{11}^{(n)} = \{p_{12}p_{23}p_{31} + p_{13}p_{32}p_{21} + p_{12}p_{21}p_{22} + p_{13}p_{31}p_{33}\}a_{n-3} + (p_{12}p_{21} + p_{13}p_{31})(p_{23}p_{32} - p_{22}p_{33})a_{n-4},$$
(4.16)
where $a_0 = 1, a_1 = p_{22} + p_{33}, a_2 = p_{22}^2 + p_{33}^2 + p_{22}p_{33} + p_{23}p_{32},$
$a_3 = (p_{22} + p_{33})(p_{22}^2 + p_{33}^2 + 2p_{23}p_{32}),$
$a_4 = (p_{22} + p_{33})^2\{(p_{22} - p_{33})^2 + p_{22}p_{33} + 3p_{23}p_{32}\} + (p_{23}p_{32} - p_{22}p_{33})^2,$
and, for $n \geq 2$, $a_n = (p_{22} + p_{33})a_{n-1} + (p_{23}p_{32} - p_{22}p_{33})a_{n-2}.$ (4.17)

**Proof:**
From Eqn. (4.15) observe that $F_{11}(s)$ can be expressed as
$$F_{11}(s) = \sum_{n=0}^{\infty} f_{11}^{(n)}s^n = 1 - \{1 + (A-1)s + (B-A)s^2 - Bs^3\}\{\sum_{n=o}^{\infty} a_n s^n\},$$

where, from Eqn.(4.14),

$$a_n = C\lambda_{12}^n + D\lambda_{13}^n = \frac{\lambda_{12}^{n+1} - \lambda_{13}^{n+1}}{\lambda_{12} - \lambda_{13}} = \frac{1 - (\lambda_{13}/\lambda_{12})^{n+1}}{1 - (\lambda_{13}/\lambda_{12})} \lambda_{12}^n = \sum_{k=0}^{n} \lambda_{13}^k \lambda_{12}^{n-k}. \qquad (4.18)$$

It is easily verified that the expressions and the recurrence relationship between $a_n$, $a_{n-1}$ and $a_{n-2}$, as given by Eqn.(4.17), follow using Eqns.(4.13).

To obtain the expressions for $f_{11}^{(n)}$ observe that from Eqns. (4.18), (4.11) and (4.13),

$f_{11}^{(0)} = \text{coeff of } s^0 = 1 - a_0 = 0,$

$f_{11}^{(1)} = \text{coeff of } s^1 = -a_1 + (1-A)a_0 = -(p_{22} + p_{33}) + (p_{11} + p_{22} + p_{33}) = p_{11},$

$f_{11}^{(2)} = \text{coeff of } s^2 = -a_2 + (1-A)a_1 + (A-B)a_0 = p_{12}p_{21} + p_{13}p_{31},$

$f_{11}^{(3)} = \text{coeff of } s^3 = -a_3 + (1-A)a_2 + (A-B)a_1 + Ba_0$

$\quad = p_{12}p_{22}p_{21} + p_{13}p_{33}p_{31} + p_{12}p_{23}p_{31} + p_{13}p_{32}p_{21}.$

In general, from Eqn. (4.18), that for $n \geq 4$,

$f_{11}^{(n)} = \text{coeff of } s^n = -a_n + (1-A)a_{n-1} + (A-B)a_{n-2} + Ba_{n-3}.$

$\quad = -a_n + (p_{11} + p_{22} + p_{33})a_{n-1} + (p_{12}p_{21} + p_{23}p_{32} + p_{13}p_{31} - p_{11}p_{22} - p_{22}p_{33} - p_{33}p_{11})a_{n-2}$

$\quad + (p_{12}p_{23}p_{31} + p_{13}p_{32}p_{21} + p_{11}p_{22}p_{33} - p_{12}p_{21}p_{33} - p_{13}p_{31}p_{22} - p_{11}p_{23}p_{23})a_{n-3}.$

Further simplification, using Eqn. (4.17), yields the expression given by Eqn. (4.16).

Note that the expressions for $f_{11}^{(1)}, f_{11}^{(2)}$ and $f_{11}^{(3)}$ as given in Theorem 4.1 and Corollary 4.2 also follow from sample path arguments. We verify that expression (4.16), when $n = 4$, also leads to an expression for the probability that the recurrence time of state 1 occurs at the fourth step. This can be derived by sample path arguments. Consider all possible paths between the sets of states on successive trials, i.e. $\{1\} \to \{2,3\} \to \{2,3\} \to \{2,3\} \to \{1\}$. Thus

$f_{11}^{(4)} = p_{12}p_{22}p_{22}p_{21} + p_{12}p_{23}p_{32}p_{21} + p_{12}p_{22}p_{23}p_{31} + p_{12}p_{23}p_{33}p_{31}$

$\quad + p_{13}p_{32}p_{22}p_{21} + p_{13}p_{33}p_{32}p_{21} + p_{13}p_{32}p_{23}p_{31} + p_{13}p_{33}p_{33}p_{33}.$

Needless to say, equivalent expressions for the recurrence time distributions for the other states occur with a permutation of the indices.

**Theorem 4.3:** (The distribution of $T_{12}$ for 3-state Markov chains).

The distribution $\{f_{12}^{(n)}\}$ where $f_{12}^{(n)} = P\{T_{12} = n\}$ is given by

$f_{12}^{(1)} = p_{12}, f_{12}^{(2)} = p_{11}p_{12} + p_{13}p_{32}, f_{12}^{(3)} = p_{11}p_{11}p_{12} + p_{11}p_{13}p_{32} + p_{13}p_{33}p_{32} + p_{13}p_{31}p_{12},$

and, in general, for $n \geq 3$, provided $\delta_2 > 0$,

$$f_{12}^{(n)} = \frac{a(\lambda_{21})}{\lambda_{21} - \lambda_{23}} \lambda_{21}^{n-1} - \frac{a(\lambda_{23})}{\lambda_{21} - \lambda_{23}} \lambda_{23}^{n-1}, \qquad (4.19)$$

where $\lambda_{21} = \dfrac{p_{11} + p_{33} + \delta_2}{2}, \lambda_{23} = \dfrac{p_{11} + p_{33} - \delta_2}{2}$, with $\delta_2 = \sqrt{(p_{11} - p_{33})^2 + 4p_{13}p_{31}}$,

and $a(\lambda) = p_{12}\lambda + p_{13}p_{32} - p_{12}p_{33}.$ $\qquad (4.20)$

**Proof:** From Eqn.(4.8), $F_{12}(s) = \dfrac{a_{12}(s)}{a_{22}(s)}$ implying, from Eqns.(4.6) and (4.7),

$$F_{12}(s) = \frac{sp_{12} + s^2(p_{13}p_{32} - p_{12}p_{33})}{1 - s(p_{11} + p_{33}) + s^2(p_{11}p_{33} - p_{13}p_{31})} = \{sp_{12} + s^2(p_{13}p_{32} - p_{12}p_{33})\}\sum_{n=0}^{\infty}(F\lambda_{21}^n + G\lambda_{23}^n)s^n$$

where $\dfrac{1}{a_{22}(s)} = \dfrac{1}{(1-\lambda_{21}s)(1-\lambda_{23}s)} = \dfrac{F}{(1-\lambda_{21}s)} + \dfrac{G}{(1-\lambda_{23}s)} = \sum_{n=0}^{\infty}(F\lambda_{21}^n + G\lambda_{23}^n)s^n$, (4.21)

leading to $\lambda_{21} + \lambda_{23} = p_{11} + p_{33}$, $\lambda_{21}\lambda_{23} = p_{11}p_{33} - p_{13}p_{31}$, $F = \dfrac{\lambda_{21}}{\lambda_{21} - \lambda_{23}}$ and $G = \dfrac{-\lambda_{23}}{\lambda_{21} - \lambda_{23}}$. (4.22)

These results imply that
$F + G = 1$,

$$F\lambda_{21} + G\lambda_{23} = \left(\frac{\lambda_{21}^2 - \lambda_{23}^2}{\lambda_{21} - \lambda_{23}}\right) = \lambda_{21} + \lambda_{23} = p_{11} + p_{33},$$

$$F\lambda_{21}^2 + G\lambda_{23}^2 = \left(\frac{\lambda_{21}^3 - \lambda_{23}^3}{\lambda_{21} - \lambda_{23}}\right) = \lambda_{21}^2 + \lambda_{23}^2 + \lambda_{21}\lambda_{23} = (\lambda_{21} + \lambda_{23})^2 - \lambda_{21}\lambda_{23} = p_{11}^2 + p_{33}^2 + p_{11}p_{33} + p_{13}p_{31}.$$

Equating the coefficients of $s^n$ for $n = 0, 1, 2$ and $3$, and using the above results we obtain
$f_{12}^{(0)} = $ coeff of $s^0 = 0$,
$f_{12}^{(1)} = $ coeff of $s^1 = p_{12}(F + G) = p_{12}$,
$f_{12}^{(2)} = $ coeff of $s^2 = p_{12}(F\lambda_{21} + G\lambda_{23}) + (p_{13}p_{32} - p_{12}p_{33})(F + G) = p_{11}p_{12} + p_{13}p_{32}$,
$f_{12}^{(3)} = $ coeff of $s^3 = p_{12}(F\lambda_{21}^2 + G\lambda_{23}^2) + (p_{13}p_{32} - p_{12}p_{33})(F\lambda_{21} + G\lambda_{23})$
$\qquad = p_{11}^2 p_{12} + p_{13}p_{31}p_{12} + p_{11}p_{13}p_{32} + p_{13}p_{33}p_{32}$,

leading to the special cases when $n = 1, 2$, and $3$.
For the general case of the Theorem when, $n \geq 3$,

$$f_{12}^{(n)} = \text{coeff of } s^n = p_{12}(F\lambda_{21}^{n-1} + G\lambda_{23}^{n-1}) + (p_{13}p_{32} - p_{12}p_{33})(F\lambda_{21}^{n-2} + G\lambda_{23}^{n-2})$$

$$= \frac{(p_{12}\lambda_{21} + p_{13}p_{32} - p_{12}p_{33})}{\lambda_{21} - \lambda_{23}}\lambda_{21}^{n-1} - \frac{(p_{12}\lambda_{23} + p_{13}p_{32} - p_{12}p_{33})}{\lambda_{21} - \lambda_{23}}\lambda_{23}^{n-1} = \frac{a(\lambda_{21})}{\lambda_{21} - \lambda_{23}}\lambda_{21}^{n-1} - \frac{a(\lambda_{23})}{\lambda_{21} - \lambda_{23}}\lambda_{23}^{n-1},$$

where $a(\lambda)$ is given by Eqn.(4.20). Further,

$$a(\lambda_{21}) = \frac{p_{12}(p_{11} - p_{33}) + 2p_{13}p_{32} + p_{12}\delta_2}{2} \text{ and } a(\lambda_{31}) = \frac{p_{12}(p_{11} - p_{33}) + 2p_{13}p_{32} - p_{12}\delta_2}{2}.$$

A recurrence relationship for $f_{12}^{(n)}$ can also be obtained, similar to that derived for $f_{11}^{(n)}$.

**Corollary 4.4:** (The distribution of $T_{12}$ for 3-state Markov chains).
The probability distribution $\{f_{12}^{(n)}\}$ where $f_{12}^{(n)} = P\{T_{12} = n\}$ is given by
$f_{12}^{(1)} = p_{12}, f_{12}^{(2)} = p_{11}p_{12} + p_{13}p_{31}$, and, for $n \geq 2$
$$f_{12}^{(n)} = p_{12}b_{n-1} + (p_{13}p_{32} - p_{12}p_{33})b_{n-2}, \tag{4.23}$$
where $b_0 = 1$, $b_1 = p_{11} + p_{33}$, $b_2 = p_{11}^2 + p_{33}^2 + p_{11}p_{33} + p_{13}p_{31}$, and, for $n \geq 2$,
$$b_n = (p_{11} + p_{33})b_{n-1} + (p_{13}p_{31} - p_{11}p_{33})b_{n-2}. \tag{4.24}$$

**Proof:**
From Eqn. (4.21) observe that $F_{12}(s)$ can be expressed as

$$F_{12}(s) = \{sp_{12} + s^2(p_{13}p_{32} - p_{12}p_{33})\}\sum_{n=o}^{\infty} b_n s^n,$$

where, $b_n = F\lambda_{21}^n + G\lambda_{23}^n = \dfrac{\lambda_{21}^{n+1} - \lambda_{23}^{n+1}}{\lambda_{21} - \lambda_{23}} = \sum_{k=0}^{n} \lambda_{21}^k \lambda_{23}^{n-k}.$

It is easily verified that the expressions and the recurrence relationship between $b_n$, $b_{n-1}$ and $b_{n-2}$, as given by Eqn.(4.23).

To obtain the expressions for $f_{12}^{(n)}$ observe that from Eqns. (4.21), (4.11) and (4.13),
$f_{12}^{(0)} = 0$, $f_{12}^{(1)} = p_{12}b_0 = p_{12}$ and
$f_{12}^{(2)} = p_{12}b_1 + (p_{13}p_{32} - p_{12}p_{33})b_0 = p_{12}(p_{11} + p_{33}) + (p_{13}p_{32} - p_{12}p_{33}) = p_{12}p_{11} + p_{13}p_{32},$
implying in general, for $n \geq 2$, the result given by Eqn.(4.23).

Interchanging states 2 and 3, from Theorem 4.3 we have an analogous expression for the distribution of $T_{13}$ We state this without a proof.

**Theorem 4.5:** (The distribution of $T_{13}$ for 3-state Markov chains).
The distribution $\{f_{13}^{(n)}\}$, where $f_{13}^{(n)} = P\{T_{13} = n\}$, is given by
$f_{13}^{(1)} = p_{13}, f_{13}^{(2)} = p_{11}p_{13} + p_{12}p_{23}$, $f_{13}^{(3)} = p_{11}p_{11}p_{13} + p_{11}p_{12}p_{23} + p_{12}p_{22}p_{23} + p_{12}p_{21}p_{13}$,
and, in general, for $n \geq 3$, provided $\delta_3 > 0$,

$$f_{13}^{(n)} = \frac{b(\lambda_{31})}{\lambda_{31} - \lambda_{32}}\lambda_{31}^{n-1} - \frac{b(\lambda_{32})}{\lambda_{31} - \lambda_{32}}\lambda_{32}^{n-1}, \tag{4.25}$$

where $\lambda_{31} = \dfrac{p_{11} + p_{22} + \delta_3}{2}, \lambda_{32} = \dfrac{p_{11} + p_{22} - \delta_3}{2}$, with $\delta_3 = \sqrt{(p_{11} - p_{22})^2 + 4p_{12}p_{21}}$,

and $b(\lambda) = p_{13}\lambda + p_{12}p_{23} - p_{13}p_{22}.$ (4.26)

We are now ready to obtain a general form of the distribution of the mixing time random variable $T_1^{(0)}$. The distributions of $T_2^{(0)}$ and $T_3^{(0)}$ will follow by similar arguments.

**Theorem 4.6:** (The distribution of $T_1^{(0)}$ for 3-state Markov chains)
The probability distribution of $T_1^{(0)}$, $f_{1,n} = P\{T_1^{(0)} = n | X_0 = 1\}$ is given by
$f_{1,0} = \pi_1, f_{1,1} = \pi_2 p_{12} + \pi_3 p_{13}, f_{1,2} = \pi_2(p_{11}p_{12} + p_{13}p_{31}) + \pi_3(p_{11}p_{13} + p_{12}p_{23}),$
$f_{1,3} = \pi_2(p_{11}p_{11}p_{12} + p_{11}p_{13}p_{32} + p_{13}p_{33}p_{32} + p_{13}p_{31}p_{12})$
$\quad + \pi_3(p_{11}p_{11}p_{13} + p_{11}p_{12}p_{23} + p_{12}p_{22}p_{23} + p_{12}p_{21}p_{13}),$
and, for $n \geq 3$, provided $\delta_2 > 0, \delta_3 > 0,$

$$f_{1,n} = \pi_2\left(\frac{a(\lambda_{21})}{\lambda_{21} - \lambda_{23}}\lambda_{21}^{n-1} - \frac{a(\lambda_{23})}{\lambda_{21} - \lambda_{23}}\lambda_{23}^{n-1}\right) + \pi_3\left(\frac{b(\lambda_{31})}{\lambda_{31} - \lambda_{32}}\lambda_{31}^{n-1} - \frac{b(\lambda_{32})}{\lambda_{31} - \lambda_{32}}\lambda_{32}^{n-1}\right), \tag{4.27}$$

where $\lambda_{21} = \dfrac{p_{11} + p_{33} + \delta_2}{2}, \lambda_{23} = \dfrac{p_{11} + p_{33} - \delta_2}{2}$, with $\delta_2 = \sqrt{(p_{11} - p_{33})^2 + 4p_{13}p_{31}}$,

$\lambda_{31} = \dfrac{p_{11} + p_{22} + \delta_3}{2}, \lambda_{32} = \dfrac{p_{11} + p_{22} - \delta_3}{2}$, with $\delta_3 = \sqrt{(p_{11} - p_{22})^2 + 4p_{12}p_{21}}$,

and $a(\lambda) = p_{12}\lambda + p_{13}p_{32} - p_{12}p_{33}$, $b(\lambda) = p_{13}\lambda + p_{12}p_{23} - p_{13}p_{22}$.

**Proof:** While we can use Eqn. (2.12) for $i = 1$ and the expressions from Eqn. (4.7)
$$f_1(s) = \pi_1 + \pi_2 \frac{a_{12}(s)}{a_{22}(s)} + \pi_3 \frac{a_{13}(s)}{a_{33}(s)},$$
since we have effectively extracted the coefficients of $s^n$ for each of these separate components, we can simply use the basic result for the distribution of the mixing times from Theorem 2.1: $f_{1,0} = \pi_1$, and for $n \geq 1$, $f_{1,n} = \pi_2 f_{12}^{(n)} + \pi_3 f_{13}^{(n)}$ and the results of Theorems 4.3 and 4.5. The results of the theorem follow, with no simplification.

The derivation of the distribution of $T_1^{(1)}$ requires the knowledge the distribution of $T_{11}$, the first return to state 1, or the recurrence time of state 1. The efforts that we have gone to in deriving the results of Theorem 4.1 can now be made use of in determining the distribution of $T_1^{(1)}$.

**Theorem 4.7:** (The distribution of $T_1^{(1)}$ for 3-state Markov chains)
The probability distribution of $T_1^{(1)}$, $g_{1,n} = P\{T_1^{(1)} = n | X_0 = 1\}$, is given by
$$g_{1,1} = \pi_1 p_{11} + \pi_2 p_{12} + \pi_3 p_{13}$$
$$g_{1,2} = \pi_1(p_{12}p_{21} + p_{13}p_{31}) + \pi_2(p_{11}p_{12} + p_{13}p_{31}) + \pi_3(p_{11}p_{13} + p_{12}p_{23}),$$
and, for $n \geq 3$, provided $\delta_1 > 0$, $\delta_2 > 0$ and $\delta_3 > 0$,
$$g_{1,n} = \pi_1 \left( \frac{c(\lambda_{12})}{\lambda_{12} - \lambda_{13}} \lambda_{12}^{n-2} - \frac{c(\lambda_{13})}{\lambda_{12} - \lambda_{13}} \lambda_{13}^{n-2} \right) + \pi_2 \left( \frac{a(\lambda_{21})}{\lambda_{21} - \lambda_{23}} \lambda_{21}^{n-1} - \frac{a(\lambda_{23})}{\lambda_{21} - \lambda_{23}} \lambda_{23}^{n-1} \right)$$
$$+ \pi_3 \left( \frac{b(\lambda_{31})}{\lambda_{31} - \lambda_{32}} \lambda_{31}^{n-1} - \frac{b(\lambda_{32})}{\lambda_{31} - \lambda_{32}} \lambda_{32}^{n-1} \right), \tag{4.28}$$
where
$$\lambda_{12} = \frac{p_{22} + p_{33} + \delta_1}{2}, \lambda_{13} = \frac{p_{22} + p_{33} - \delta_1}{2}, \text{ with } \delta_1 = \sqrt{(p_{22} - p_{33})^2 + 4p_{23}p_{32}},$$
$$\lambda_{21} = \frac{p_{11} + p_{33} + \delta_2}{2}, \lambda_{23} = \frac{p_{11} + p_{33} - \delta_2}{2}, \text{ with } \delta_2 = \sqrt{(p_{11} - p_{33})^2 + 4p_{13}p_{31}},$$
$$\lambda_{31} = \frac{p_{11} + p_{22} + \delta_3}{2}, \lambda_{32} = \frac{p_{11} + p_{22} - \delta_3}{2}, \text{ with } \delta_3 = \sqrt{(p_{11} - p_{22})^2 + 4p_{12}p_{21}},$$
$$a(\lambda) = p_{12}\lambda + p_{13}p_{32} - p_{12}p_{33}, \quad b(\lambda) = p_{13}\lambda + p_{12}p_{23} - p_{13}p_{22},$$
$$c(\lambda) = (1 - p_{11})(1 - p_{22})(1 - p_{33}) - p_{23}p_{32}(1 - p_{11}) - (p_{12}p_{21} + p_{13}p_{31})(1 - \lambda).$$

**Proof:** From Eqn.(2.2), for $n \geq 1$, $g_{1,n} = \pi_1 f_{11}^{(n)} + \pi_2 f_{12}^{(n)} + \pi_3 f_{13}^{(n)}$ with $g_{1,0} = 0$.
Using the results of Theorems 4.1, 4.3 and 4.5 the results of the theorem follow directly.

Before we examine some special cases, we have the following results that we state without proof, The results are given in [3] and [5], although alternative proofs can be given using the results of this paper (analogous to the proof of Theorem 3.3 above).

**Theorem 4.8:** (Mean mixing times for 3-state Markov chains)
If the three-state Markov chain with transition matrix given by Eqn. (4.1) is irreducible, (if and only if $\Delta_1 > 0$, $\Delta_2 > 0$, and $\Delta_3 > 0$), then

$$E\left[T_1^{(0)}\right] = E\left[T_2^{(0)}\right] = \frac{\kappa}{\Delta} = \tau \tag{4.28}$$

$$E\left[T_1^{(1)}\right] = E\left[T_2^{(1)}\right] = 1 + \frac{\kappa}{\Delta} = \eta \tag{4.29}$$

where $\kappa = p_{12} + p_{13} + p_{21} + p_{23} + p_{31} + p_{32}$, (4.30)

and $\Delta = p_{23}p_{31} + p_{21}p_{32} + p_{21}p_{31} + p_{31}p_{12} + p_{32}p_{13} + p_{32}p_{12} + p_{12}p_{23} + p_{13}p_{21} + p_{13}p_{23}$. (4.31)

To illustrate the general results of Theorems 4.6 and 4.7 we consider some special cases of three state Markov chains. These cases were considered in [5] to illustrate some general results for the expected values of the mixing time random variable $T_i^{(1)}$ and in [3] when the expected value of the alternative mixing random variable $T_i^{(0)}$ was considered.

*Case* 1: *"Minimal period 3"*

Let $P = \begin{bmatrix} 0 & 1 & 0 \\ 0 & 0 & 1 \\ 1 & 0 & 0 \end{bmatrix}$, implying that the Markov chain is periodic, period 3, with transitions occurring $1 \to 2 \to 3 \to 1$ .... Then $\Delta_1 = \Delta_2 = \Delta_3 = 1$, $\Delta = 3$, and $\pi_1 = \pi_2 = \pi_3 = 1/3$.
It is easily seen that $f_{11}^{(3)} = 1, f_{12}^{(2)} = 1, f_{13}^{(2)} = 1$, leading to

$$f_{1,n} = \frac{1}{3}, n = 0,1,2,$$

and
$$g_{1,n} = \frac{1}{3}, n = 1,2,3.$$

This is consistent with the following observations. Suppose that the Markov chain starts in state 1, with $X_0 = 1$. If the mixing state $M$ is 1 (with probability 1/3) then $T_1^{(0)} = 0$, so that mixing occurs at that trial while $T_1^{(1)} = 3$ since 3 further steps $1 \to 2$, $2 \to 3$, $3 \to 1$. If the mixing state $M$ is 2 (with probability 1/3) then $T_1^{(0)} = T_1^{(1)} = 1$ since the mixing state occurs after 1 further step as $1 \to 2$. If the mixing state $M$ is 3 (with probability 1/3) then $T_1^{(0)} = T_1^{(1)} = 2$ since the mixing state occurs after 2 further steps as $1 \to 2$, $2 \to 3$.

A simple deduction is that $E[T_1^{(0)}] = 1$ and $E[T_1^{(1)}] = 2$, consistent with the observations of Theorem 2.7 and the earlier result, reported in [4,] that for irreducible periodic, period *3*, Markov chains, $\eta = 2$, being the minimal value of the expected time to mixing in a three state Markov chain.

*Case 2: "Period 2"*

Let $P = \begin{bmatrix} 0 & 1 & 0 \\ q & 0 & p \\ 0 & 1 & 0 \end{bmatrix}$, ($p + q = 1$), the transition matrix of a periodic period 2 three-state Markov chain (with transitions alternating between the states {1, 3} and {2}).

Then $\Delta_1 = q$, $\Delta_2 = 1$, $\Delta_3 = p$, $\Delta = 2$ implying $\pi_1 = q/2$, $\pi_2 = 1/2$, $\pi_3 = p/2$.

In Example 6.1.7 of [7] (whilst obtaining explicit expressions for the *n*-step transition probabilities) an expressions for $\mathbf{P}(s)$ was shown to be

$$\mathbf{P}(s) = \frac{1}{(1-s^2)} \begin{bmatrix} 1-ps^2 & s & ps^2 \\ qs & 1 & ps \\ qs^2 & s & 1-qs^2 \end{bmatrix}.$$

This leads immediately to the results (using Eqn. (4.8)) that

$$F_{11}(s) = \frac{qs^2}{1-ps^2}, F_{12}(s) = s, F_{13}(s) = \frac{ps^2}{1-qs^2}.$$

Extraction of the coefficients of $s^n$ (via power series expansions for $F_{11}(s)$ and $F_{13}(s)$), lead to $f_{12}^{(1)} = 1$ (consistent with the observation that a path from 1 always leads in one step to 2), and $f_{11}^{(2n+2)} = qp^n, f_{13}^{(2n+2)} = pq^n, (n \geq 0)$. This implies that the mixing time distributions are given as

$$f_{1,0} = \frac{q}{2}, f_{1,1} = \frac{1}{2}, f_{1,2n+2} = \frac{p^2 q^n}{2} \ (n \geq 0),$$

and

$$g_{1,1} = \frac{1}{2}, g_{1,2n+2} = \frac{q^2 p^n + p^2 q^n}{2}, (n \geq 0).$$

*Case 3: "Constant movement"*

Let $P = \begin{bmatrix} 0 & b & c \\ d & 0 & f \\ g & h & 0 \end{bmatrix}$, ($b+c=1$, $d+f=1$, $g+h=1$). In this case $p_{11} = p_{22} = p_{33} = 0$, so that at each step the chain does not remain at the state but moves to one of the other states. The Markov chain is irreducible, and regular if $0 < b < 1, 0 < f < 1, 0 < g < 1$.

Now $\Delta_1 = 1 - fh$, $\Delta_2 = 1 - cg$, $\Delta_3 = 1 - bd$, $\Delta = 3 - fg - gc - bd$,

implying $\pi_1 = \frac{1-fh}{3-fh-cg-bd}$, $\pi_2 = \frac{1-cg}{3-fh-cg-bd}$, and $\pi_3 = \frac{1-bd}{3-fh-cg-bd}$.

In [5] it was shown that $1 \leq E[T_i^{(1)}] = \eta \leq 1.5$. The minimal value of $\eta = 1$ occurs when either $b = f = g = 1$ (and this case reduces to the "period 3" Case 1 above), or when $b = f = g = 0$ (when this case again reduces to a periodic, "period 3" chain but with transitions $1 \to 3 \to 2 \to 1$ ....).

The maximal value of $\eta = 1.5$ occurs when any pair of ($b, f, g$) take the values 0 and 1, say $b = 1, g = 0$, when this case reduces to the "period 2" Case 2 above.

For the regular case $1 < \eta < 1.5$, which we now explore.

After simplification of the algebra, from Theorem 4.1, the distribution $\{f_{11}^{(n)}\}$ is given by

$f_{11}^{(1)} = 0; f_{11}^{(2n)} = (bd + cg)(fh)^{n-1}, (n = 1, 2, ...); f_{11}^{(2n+1)} = (1 - fh - bd - cg)(fh)^{n-1}, (n = 1, 2, ...).$

From Theorem 4.3, the distribution $\{f_{12}^{(n)}\}$ is given by

$f_{12}^{(2n)} = ch(cg)^{n-1}, (n = 1, 2, ...); f_{12}^{(2n+1)} = b(cg)^n, (n = 0, 1, 2, ...).$

From Theorem 4.5, the distribution $\{f_{13}^{(n)}\}$ is given by

$f_{13}^{(2n)} = bf(bd)^{n-1}, (n = 1, 2, ...); f_{13}^{(2n+1)} = c(bd)^n, (n = 0, 1, 2, ...).$

From Theorem 2.1, the distribution of the mixing time random variable $T_1^{(0)}, \{f_{1,n}\}$, is given by

$$f_{1,0} = \frac{1 - fh}{3 - fh - cg - bd}, \quad f_{1,1} = \frac{1 - bc(d + g)}{3 - fh - cg - bd}$$

$$f_{1,2n} = \frac{ch(1 - cg)}{3 - fh - cg - bd}(cg)^{n-1} + \frac{bf(1 - bd)}{3 - fh - cg - bd}(bd)^{n-1}, \ (n = 1, 2, ...),$$

$$f_{1,2n+1} = \frac{b(1 - cg)}{3 - fh - cg - bd}(cg)^n + \frac{c(1 - bd)}{3 - fh - cg - bd}(bd)^n, \ (n = 1, 2, ...) \text{ and also } n = 0, \text{ as above,}$$

The distribution of the mixing time random variable $T_1^{(1)}, \{g_{1,n}\}$, is given by

$$g_{1,0} = 0, \ g_{1,1} = \frac{1 - bc(d + g)}{3 - fh - cg - bd},$$

$$g_{1,2n} = \frac{(bd + cg)(1 - fh)(fh)^{n-1}}{3 - fh - cg - bd} + \frac{ch(1 - cg)(cg)^{n-1}}{3 - fh - cg - bd} + \frac{bf(1 - bd)(bd)^{n-1}}{3 - fh - cg - bd}, (n = 1, 2, ...)$$

$$g_{1,2n+1} = \frac{(1 - fh - bd - cg)(1 - fh)(fh)^{n-1}}{3 - fh - cg - bd} + \frac{b(1 - cg)(cg)^{n-1}}{3 - fh - cg - bd} + \frac{c(1 - bd)(bd)^{n-1}}{3 - fh - cg - bd}, (n = 1, 2, ...)$$

Thus the mixing time distributions are basically mixtures of modified geometric distributions.

We consider the special case of $b = f = g = \varepsilon, c = d = h = 1 - \varepsilon$.

$$f_{1,0} = \frac{1}{3}, \ f_{1,2n} = \frac{\{(1-\varepsilon)^2 + \varepsilon^2\}\{\varepsilon(1-\varepsilon)\}^{n-1}}{3}, \ (n = 1, 2, ...),$$

$$f_{1,2n+1} = \frac{\{\varepsilon(1-\varepsilon)\}^n}{3}, \ (n = 0, 1, 2, ...);$$

and

$$g_{1,0} = 0, \ g_{11} = \frac{1}{3}, g_{1,2n} = \frac{\{\varepsilon(1-\varepsilon)\}^{n-1}}{3}, (n = 1, 2, ...),$$

$$g_{1,2n+1} = \frac{\{(1-\varepsilon)^2 + \varepsilon^2\}\{\varepsilon(1-\varepsilon)\}^{n-1}}{3}, (n = 1, 2, ...)$$

It is interesting to observe, in comparing the two mixing time distributions $\{f_{1,n}\}$ and $\{g_{1,n}\}$, that, for $n = 0, 1, 2, \ldots$ $f_{1,2n} = g_{1,2n+1}$ and $f_{1,2n+1} = g_{1,2n+2}$, so that for all $n \geq 0$,
$P\{T_1^{(0)} = n\} = f_{1,n} = g_{1,n+1} = P\{T_1^{(1)} = n+1\} = P\{T_1^{(1)} - 1 = n\}$. Thus $T_1^{(0)}$ has the same distribution as $T_1^{(1)} - 1$, or equivalently that $T_1^{(1)}$ is distributed as $T_1^{(0)} + 1$. While we have earlier shown that in every mixing situation, $E[T_1^{(1)}] = E[T_1^{(0)}] + 1$, this does not necessarily imply, in general, that $T_1^{(1)}$ and $T_1^{(0)} + 1$ have the same distribution, as observed for the above situation.

*Case 4: "Independent"*

Let $P = \begin{bmatrix} a_1 & a_2 & a_3 \\ a_1 & a_2 & a_3 \\ a_1 & a_2 & a_3 \end{bmatrix}$, $(a_1 + a_2 + a_3 = 1)$, implying that the Markov chain is equivalent to independent trials on the state space $S = \{1, 2, 3\}$.

Observe that $\Delta_1 = a_1$, $\Delta_2 = a_2$, $\Delta_3 = a_3$, $\Delta = 1$ implying $\pi_1 = a_1$, $\pi_2 = a_2$, $\pi_3 = a_3$.
It is easily seen that for $j = 1, 2, 3$, $f_{1j}^{(n)} = a_j(1-a_j)^{n-1}$, $n = 1,2,3,\ldots$ (i.e. geometric $(a_j)$ distributions) implying that the two mixing time distributions are given by
$f_{10} = a_1, f_{1,n} = a_2^2(1-a_2)^{n-1} + a_3^2(1-a_3)^{n-1}$, $n = 2,3,\ldots$,
and $g_{10} = 0, g_{1,n} = a_1^2(1-a_1)^{n-1} + a_2^2(1-a_2)^{n-1} + a_3^2(1-a_3)^{n-1}$, $n = 1,2,3,\ldots$,
i.e. a mixture of three distributions – a constant and two geometric distributions for $\{f_{1,n}\}$, and three geometric distributions for $\{g_{1,n}\}$.

In Case 3 the mixing random variables have relatively tight distributions since the mean times to mixing are constrained within tight bounds. We finish with a case where the mixing time random variables can take relatively large values, by constraining the movement within the states to ensure that the Markov chains can reside in individual states for possibly long periods of time before moving.

*Case 5: "Cyclic drift"*

Let $P = \begin{bmatrix} a & b & 0 \\ 0 & c & d \\ f & 0 & g \end{bmatrix}$, $a+b=1$, $c+d=1$, $f+g=1$, with $0 < a < 1$, $0 < c < 1$, $0 < g < 1$,

implying that the Markov chain is regular. Observe that at each transition the chain either remains in the same state $i$ or moves to state $i + 1$ (or 1 if $i = 3$).
Now $\Delta_1 = df$, $\Delta_2 = fb$, $\Delta_3 = bd$ so that
$\pi_1 = \dfrac{df}{df + fb + bd}$, $\pi_2 = \dfrac{fb}{df + fb + bd}$, and $\pi_3 = \dfrac{bd}{df + fb + bd}$.

From [3] and [5], $E[T_1^{(1)}] = \eta = 1 + \frac{b+d+f}{df+bf+bd}$, $E[T_1^{(0)}] = \tau = \frac{b+d+f}{df+bf+bd}$.

Note that $0 < b+d+f < 3$ and $0 < df+bf+bd < 3$.

When $b+d+f \to 3$ then $df+bf+bd \to 3$ and $\eta \to 2$ (as in Case 1).

When $b+d+f \to 0$ then $df+bf+bd \to 0$ but the behaviour of $\eta$ and $\tau$ depends upon the rates of convergence.

Let $b = d = f = \varepsilon$, then $2 < \eta = 1 + \frac{1}{\varepsilon} < \infty$ where the lower bound is achieved in the periodic, non-regular case ($\varepsilon = 1$), as in Case 1. Arbitrary large values of $\eta$ occur as $\varepsilon \to 0$ when the Markov chain is approaching the reducible situation with all states absorbing.

What is emerging is that if the Markov chain has states where it resides for a large number of transitions, i.e. if there is little movement the mixing time can become excessively large. We explore this in more detail.

Firstly, we evaluate the recurrence time and first passage time distributions from state 1, using Theorems 4.1, 4.3 and 4.5.

Provided $c \neq g$, the distribution of $T_{11}, \{f_{11}^{(n)}\}$ is given by

$f_{11}^{(1)} = a; f_{11}^{(2)} = 0; f_{11}^{(n)} = \frac{bdf(c^{n-2} - g^{n-2})}{c-g}, n = 3, 4, \dots$.

When $c = g$ (and thus $d = 1 - c = 1 - g = f$)

$f_{11}^{(1)} = a; f_{11}^{(2)} = 0; f_{11}^{(n)} = (n-2)b(1-c)^2 c^{n-3}, n = 3, 4, \dots$.

From Theorem 4.3, or sample path arguments, the distribution of $T_{12}, \{f_{12}^{(n)}\}$ is given by

$f_{12}^{(n)} = ba^{n-1}, (n = 1, 2, \dots)$.

From Theorem 4.5, provided $a \neq c$ the distribution of $T_{13}, \{f_{13}^{(n)}\}$ is given by

$f_{13}^{(1)} = 0; f_{13}^{(2)} = bd; f_{13}^{(n)} = \frac{bd(a^{n-1} - c^{n-1})}{a-c}, n = 3, 4, \dots$.

When $a = c$ (and thus $b = 1 - a = 1 - c = d$)

$f_{13}^{(1)} = 0, f_{13}^{(n)} = (n-1)(1-a)^2 a^{n-2}, n = 2, 3, \dots$.

Let us assume that $a \neq c \neq g$. Then the two mixing time distributions are

$f_{1,0} = \frac{df}{df + fb + bd}, f_{1,1} = \frac{fb^2}{df + fb + bd}, f_{1,2} = \frac{b^2(af + d^2)}{df + fb + bd}$,

$f_{1,n} = \frac{fb^2 a^{n-1}}{df + fb + bd} + \frac{b^2 d^2 (a^{n-1} - c^{n-1})}{(df + fb + bd)(a-c)}, n = 3, 4, \dots$.

and

$g_{1,0} = 0, g_{1,1} = \frac{adf + b^2 f}{df + fb + bd}, g_{1,2} = \frac{ab^2 f + b^2 d^2}{df + fb + bd}$,

$g_{1,n} = \frac{bd^2 f^2 (c^{n-2} - g^{n-2})}{(df + fb + bd)(c-g)} + \frac{fb^2 a^{n-1}}{df + fb + bd} + \frac{b^2 d^2 (a^{n-1} - c^{n-1})}{(df + fb + bd)(a-c)}, n = 3, 4, \dots$.

Of interest is the special case when $b = d = f = \varepsilon$, $a = c = g = 1 - \varepsilon$, implying that $\pi_1 = \pi_2 = \pi_3 = \frac{1}{3}$.

In this case

$f_{11}^{(1)} = 1 - \varepsilon; f_{11}^{(2)} = 0; f_{11}^{(n)} = (n-2)\varepsilon^3(1-\varepsilon)^{n-3}, n = 3, 4, \ldots$ with $m_{11} = E[T_{11}] = 3$;

$f_{12}^{(n)} = \varepsilon(1-\varepsilon)^{n-1}, (n = 1, 2, \ldots)$ with $m_{12} = E[T_{12}] = \frac{1}{\varepsilon}$;

$f_{13}^{(1)} = 0, f_{13}^{(n)} = (n-1)\varepsilon^2(1-\varepsilon)^{n-2}, n = 2, 3, \ldots$, with $m_{13} = E[T_{13}] = \frac{2}{\varepsilon}$.

The mixing time distributions are given by

$f_{1,0} = \frac{1}{3}, f_{1,1} = \frac{\varepsilon}{3}, f_{1,2} = \frac{\varepsilon}{3}, f_{1,3} = \frac{\varepsilon(1-\varepsilon^2)}{3}, f_{1,n} = \frac{\varepsilon(1-\varepsilon)^{n-2}\{1+(n-2)\varepsilon\}}{3}, n = 4, \ldots$;

$g_{1,0} = 0, g_{1,1} = \frac{1}{3}, g_{1,2} = \frac{\varepsilon}{3}, g_{1,3} = \frac{\varepsilon}{3}, g_{1,n} = \frac{\varepsilon(1-\varepsilon)^{n-3}\{1+(n-3)\varepsilon\}}{3}, n = 4, \ldots$.

The expected mixing times are $E[T_1^{(1)}] = \eta = 1 + \frac{1}{\varepsilon}$ and $E[T_1^{(0)}] = \tau = \frac{1}{\varepsilon}$.

Graph 1 gives a plot of nine different variants of $f_{1,n}$ for $n = 0(1)20$ and as $\varepsilon$ takes the nine values $0.1(0.1)0.9$. We see that probability distribution places increasing weight on the "tail probabilities" as ε decreases. This is a reflection of the increasing mean of the distribution as ε decreases, since $E[T_1^{(0)}] = \sum_{n=1}^{\infty} nf_{1,n} = \frac{1}{\varepsilon}$.

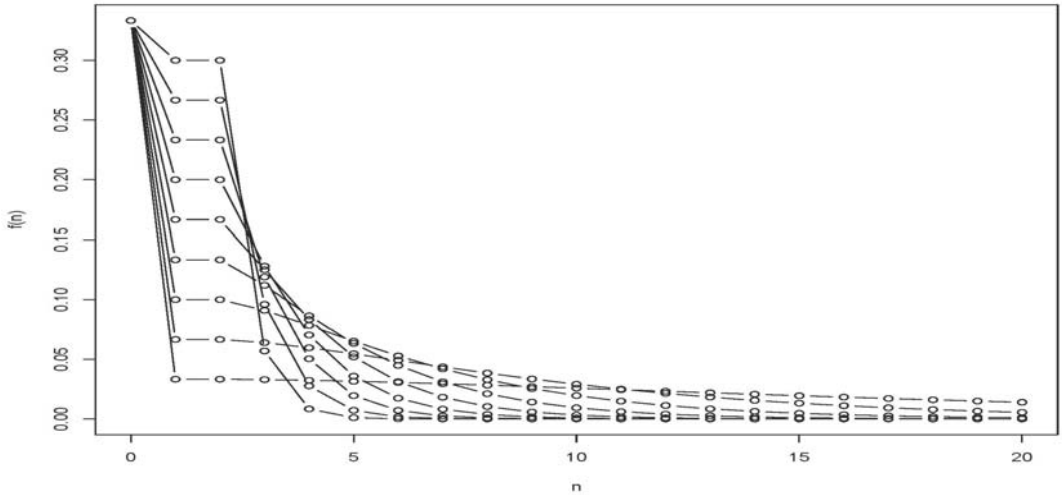

*Graph 1: Plot of the mixing time distribution {$f_{1,n}$}*

The observation that we made in comparing the two mixing time distributions $\{f_{1,n}\}$ and $\{g_{1,n}\}$ in the special case of Case 3 also holds here since for $n = 0, 1, 2, \ldots$ $P\{T_1^{(0)} = n\} = f_{1,n} = g_{1,n+1} = P\{T_1^{(1)} = n+1\}$ implying that $T_1^{(1)}$ is distributed as $T_1^{(0)} + 1$.